\documentclass[article,oneside,reqno,british,pdftex]{memoir}

\usepackage{color}
\definecolor{darkblue}{rgb}{0,0.08,0.50} 
\usepackage{amsmath,amsthm,amssymb}
\usepackage[british]{babel} 
\usepackage{microtype} 
\usepackage[modulo,mathlines]{lineno}
\usepackage[sort,square]{natbib} 
\usepackage[bookmarksnumbered=true,colorlinks=true, linkcolor=darkblue,
            citecolor=darkblue, urlcolor=darkblue, pdfpagemode=none]
           {hyperref} 
\usepackage{memhfixc} 
\DeclareFontFamily{OT1}{rsfs}{}
\DeclareFontShape{OT1}{rsfs}{n}{it}{<-> rsfs10}{}
\DeclareMathAlphabet{\curly}{OT1}{rsfs}{n}{it}
\newcommand{\smallspace}{\vspace{2mm}}

\newcommand{\reg}{\operatorname{reg}}

\newcommand\PP{\mathbb P}

\renewcommand\O{\mathcal O}
\newcommand\I{\curly I}

\newcommand{\linkref}[2]{\hyperref[#2]{#1 \ref{#2}}}

\newcommand{\num}{\text{num}}

\makeatletter \@addtoreset{equation}{chapter} \makeatother

\newtheorem{thm}[equation]{Theorem}
\newtheorem{lem}[equation]{Lemma}

\newtheorem{cor}[equation]{Corollary}
\newtheorem{cor*}{Corollary}
\newtheorem{prop}[equation]{Proposition}

\theoremstyle{definition}
\newtheorem{example}[equation]{Example}
\newtheorem*{example*}{Example}
\newtheorem{question}[equation]{Question}
\newtheorem*{question*}{Question}
\newtheorem{defn}[equation]{Definition}
\newtheorem*{defn*}{Definition}
\newtheorem{rmk}[equation]{Remark}

\newtheorem*{rmk*}{Remark}
\newtheorem*{rmks*}{Remarks}

\newtheorem*{note*}{Note}

\makeatletter
\@ifpackageloaded{lineno}%
{
\newcommand*\patchAmsMathEnvironmentForLineno[1]{%
  \expandafter\let\csname old#1\expandafter\endcsname\csname #1\endcsname
  \expandafter\let\csname oldend#1\expandafter\endcsname\csname end#1\endcsname
  \renewenvironment{#1}%
     {\linenomath\csname old#1\endcsname}%
     {\csname oldend#1\endcsname\endlinenomath}}%
\newcommand*\patchBothAmsMathEnvironmentsForLineno[1]{%
  \patchAmsMathEnvironmentForLineno{#1}%
  \patchAmsMathEnvironmentForLineno{#1*}}%
\AtBeginDocument{%
\patchBothAmsMathEnvironmentsForLineno{equation}%
\patchBothAmsMathEnvironmentsForLineno{eqnarray}%
\patchBothAmsMathEnvironmentsForLineno{align}%
\patchBothAmsMathEnvironmentsForLineno{flalign}%
\patchBothAmsMathEnvironmentsForLineno{alignat}%
\patchBothAmsMathEnvironmentsForLineno{gather}%
\patchBothAmsMathEnvironmentsForLineno{multline}%
}
}%
\makeatother

\begin{document}

 \pagestyle{myheadings}
\title{Slope Stability and Exceptional Divisors of High Genus}
\author{Dmitri Panov and Julius Ross\thanks{First author supported by EPSRC grant number EP/E044859/1.  Second author partially supported by the National Science Foundation, Grant No.\  DMS-0700419}}
\date{22nd October 2007}
\maketitle

\hypersetup{
pdfauthor = {\theauthor},
pdftitle = {\thetitle},
pdfcreator = {LaTeX with hyperref package},
}


\bibliographystyle{ross} 


\begin{abstract}
  We study slope stability of smooth surfaces and its connection with exceptional divisors.
  We show that a surface containing an exceptional divisor with arithmetic genus at least two is slope unstable for some polarisation.   In the converse direction we show  that slope stability of surfaces can be tested with divisors, and prove that for surfaces with non-negative Kodaira dimension any destabilising divisor must have negative self-intersection and arithmetic genus at least two.  We also prove that a destabilising divisor can never be nef, and as an application give an example of a surface that is slope stable but not K-stable.
\end{abstract}

\pagestyle{myheadings} 
\markboth{\thetitle}{\thetitle}

\ifdraftdoc\linenumbers\else\fi 


\chapter{Introduction}
There are many notions of stability for a projective variety.  A
classical one is Chow stability that requires  the Chow point of
the variety be stable in the sense of Geometric Invariant Theory.
Another notion is
K-stability which is closely connected to the existence of
extremal K\"ahler metrics.    For higher dimensional varieties the problem of giving an intrinsic characterisation of either of these stability notions remains largely open.  

Related to K-stability is another notion called slope
stability introduced by Thomas and the second author
\cite{ross_thomas:06:obstr_to_exist_const_scalar,ross_thomas:07:study_hilber_mumfor_criter_for}.
This involves a definition of slope $\mu(X)$ for a pair
$(X,L)$ consisting of a variety $X$ and ample line bundle $L$, and a slope $\mu(\O_Z)$ for each subscheme $Z\subset
X$ (precise definitions appear below).  We say $(X,L)$ is 
\emph{slope semistable} if $\mu(X)\le \mu(\O_Z)$ for every $Z$, and if
this fails we say $Z$ \emph{destabilises} $X$.  These definitions are
made so that K-stability implies slope stability, and thus gives a concrete and geometric obstruction to K-stability.    In turn, through the work of several authors, it also gives an obstruction to the existence of constant scalar curvature K\"ahler metrics (the result we are referring to was first proved by Donaldson \cite{donaldson(01):scalar_curvat_projec_embed}, but work of Tian \cite{tian97:kaehler_einstein_metrics_positive}, Mabuchi \cite{mabuchi:05:energ_theor_approac_to_hitch} and Chen-Tian \cite{chen-tian:geomet_kaehl_metric_holom_foliat_by_discs} contain similar results in this direction; the simplest proof can be found in \cite{donaldson:05:lower_bound_calab_funct}).

A natural question arises as to which varieties are slope stable,
and if there are any restrictions on the geometry of destabilising
subschemes.    In this paper we explore a connection between slope stability of smooth surfaces and the existence of exceptional divisors of high genus (recall that an effective divisor $D=\sum d_i  D_i$ in a surface is said to be \emph{exceptional} if its intersection matrix $(D_i.D_j)$ is negative definite).

\begin{thm}\label{thm:introthm1}\ 
  Let $X$ be a smooth projective surface containing an exceptional divisor that has arithmetic genus at least two.   Then  $(X,L)$ is slope unstable for some polarisation $L$. \end{thm}

This theorem is a natural generalisation of examples of unstable surfaces due to the second author \cite{ross:06:unstab_produc_smoot_curves}, and due to Shu \cite{shu:unstab_kodair_fibrat}.  In fact the polarisations used in this Theorem make the volume of the exceptional divisor very small, and one should interpret this as saying that the space obtained by contracting this divisor is unstable due to the existence of a singularity that is not rational or elliptic (see \linkref{Remark}{rmk:singular}).\smallspace 

As a first step in the converse direction we show that for smooth surfaces slope stability can be tested with divisors.

\begin{thm}\label{thm:introthm1.1}\ 
If a smooth polarised surface $(X,L)$ is slope unstable then it is destabilised by a divisor. 
\end{thm}

Given this, the interest lies in the geometry of destabilising divisors.  When $X$ has non-negative Kodaira dimension it is easy to prove that any destabilising divisor has arithmetic genus at least two (\linkref{Proposition}{prop:genus}).  The bulk of this paper is devoted to proving the next theorem that contains two more results in this direction.

\begin{thm}\label{thm:introthm2}
Let $D$ be an effective divisor in a smooth surface $X$ and either
  \begin{enumerate}
  \item Suppose $X$ has non-negative Kodaira dimension and $D^2\ge 0$, or
  \item Suppose $D$ is nef (i.e.\  $D.C\ge 0$ for every irreducible curve $C$ in $X$).
  \end{enumerate}
Then $D$ does not destabilise $X$ with respect to any polarisation.
\end{thm}

Although there is clearly some overlap between the two statements of this theorem  the proofs we offer are different (interestingly, positivity of the ``adjoint divisor'' $K_X+nD$ plays a central role in both).   As an application we show that slope stability is not equivalent to K-stability (in fact $\PP^2$ blown up at two points witnesses this difference; see \linkref{Example}{ex:p2blown2points}). \smallspace

Any exceptional divisor that is sufficiently close to a destabilising divisor will also be destabilising.  Thus it is unreasonable to expect that a destabilising divisor is necessarily exceptional (see \eqref{ex:nonexcep} for a specific example).   However motivated by the examples we have studied we ask the following question:

\begin{question}
  Suppose that $D$ is a divisor that destabilises a polarised
  surface $(X,L)$.  Is it necessarily true that there is an exceptional divisor $D'\subset D$ that also destabilises $(X,L)$?  
\end{question}

The paper is organised as follows.  We start with a brief introduction to slope stability in \linkref{Section}{sec:slopestability}, and prove \linkref{Theorem}{thm:introthm1} in \linkref{Section}{sec:exceptional}.   The reduction of slope stability of surfaces to divisors is made in \linkref{Section}{sec:reduction}.    The proof of the statement in \linkref{Theorem}{thm:introthm2} concerning surfaces of non-negative Kodaira dimension is started in  \linkref{Section}{sec:nonnegativeI} where it is proved for minimal surfaces and completed in \linkref{Section}{sec:nonnegativeII}.  We then give a self-contained proof of the statement about nef divisors in \linkref{Section}{sec:nef}.

\subsection{Terminology}
We work throughout over the complex numbers.  We shall write
$D\equiv E$ to mean that $D$ and $E$ are numerically equivalent divisors but, when it
cannot cause confusion, we will denote a divisor and the numerical
class it represents by the same letter and simply write $D=E$.  The canonical class of $X$
will be denoted by $K_X$, and the N\'eron-Severi space $N^1(X)$ of a smooth variety $X$
is the space of divisors in $X$ modulo numerical equivalence.  

Inside $N^1(X)\otimes \mathbb Q$ lies the \emph{ample cone} consisting of numerical classes of ample $\mathbb Q$-divisors, and the \emph{effective cone} consisting of the numerical
classes of effective $\mathbb Q$-divisors.  By the Kleiman criterion, the closure
of the ample cone is the \emph{nef cone} and consists of classes that
have non-negative intersection with every irreducible curve in $X$.  Thus the nef cone is dual to the closure of the effective cone which is called the \emph{pseudoeffective} cone.
By a \emph{polarised} variety $(X,L)$ we mean a variety $X$ and choice
of ample $\mathbb Q$-divisor $L$, and refer to $L$ as the
\emph{polarisation}.

When $X$ is a surface we denote the intersection of two divisors $D$
and $E$ in $N^1(X)$ by $D.E$.  The arithmetic genus of an effective divisor
$D$ is defined to be $p_a(D) = 1-\chi(\O_D)$ and by the adjunction formula is
given by $K_X.D+D^2 = 2p_a(D)-2$.  If $D$ is an irreducible rational curve with self-intersection $-n$ we say $D$ is a $-n$ curve.

\subsection{Acknowledgements.} We would particularly like to thank Simon Donaldson 
and Alessio Corti for useful discussions and for answering many
questions about surfaces, and Richard Thomas and G\'abor Sz\'ekelyhidi
for helpful conversations and for comments on a first version of this paper.  We also thank Brian Conrad,  Johan de Jong, J\'anos Koll\'ar,  Citriona Mclean and Jacopo Stoppa.

\chapter{Slope stability}\label{sec:slopestability}

This section contains a brief summary of slope stability and the
reader is referred to
\cite{ross_thomas:06:obstr_to_exist_const_scalar,ross_thomas:07:study_hilber_mumfor_criter_for} for
more details.  

Let $X$ be a smooth variety of dimension $n$ and $L$ an ample
divisor on $X$ with Hilbert polynomial
\[\chi(\O(kL)) = a_0 k^n + a_1k^{n-1} + O(k^{n-2}).\]
The \emph{slope of $(X,L)$} is defined to be the ratio
\[\mu(X,L) = \frac{a_1}{a_0}.\]
Now suppose $Z\subset X$ is a proper closed subscheme and 
$\pi\colon Bl_Z\to X$ is the blowup along $Z$ with exceptional divisor
$E$.  For sufficiently small $c\in \mathbb Q^+$ the divisor $\pi^*L -cE$ is ample, and
the \emph{Seshadri constant} of $Z$ is 
\[\varepsilon(Z,L) = \text{sup}\{c\in \mathbb Q : \pi^* L -cE \text{ is
  ample}\}.\] Using this we define the slope of $Z$.  For fixed $x\in \mathbb Q$ define $a_i(x)$ by
\[\chi(\O(k(\pi^*L-xE))) = a_0(x) k^n + a_1(x) k^{n-1} + O(k^{n-2})
\quad\text{ for all } k\gg 0, kx\in \mathbb N.\] From standard theory
$a_i(x)$ is a polynomial in $x$ which we extend by continuity to
$x\in\mathbb R$.  Set $\tilde{a}_i(x)=a_i-a_i(x)$.  The \emph{slope of $Z$} with respect to a given real number $c$ is defined to be
\[\mu_c(\O_Z,L) = \frac{\int_0^c \tilde{a}_1(x) +
  \frac{\tilde{a}_0'(x)}{2}dx }{\int_0^c \tilde{a}_0(x) dx}.\] If $0<c\le \varepsilon(Z,L)$ then this
quantity is finite; in fact $\tilde{a}_0(0)=0$ and $\tilde{a}_0(x)$ is
increasing in the range $0<x<\varepsilon(Z)$ so the denominator is
strictly positive \cite[Rmk.\ 3.10]{ross_thomas:06:obstr_to_exist_const_scalar}.
When the polarisation is clear from context we suppress the $L$ in
the notation and write $\mu_c(\O_Z)$, $\mu(X)$ and $\varepsilon(Z)$.

\begin{defn}  We say a polarised variety $(X,L)$ is
  \emph{slope semistable with respect to $Z$} if
  \begin{equation*}
    \mu(X,L)\le \mu_c(\O_Z,L) \quad \text{ for all } 0<c\le \varepsilon(Z,L).\label{eqn:slopstablewrtZ}   
  \end{equation*}
  We say $(X,L)$ is \emph{slope semistable} if it is slope semistable with
  respect to all subschemes; when this fails  we say $Z$ \emph{destabilises} $(X,L)$.
\end{defn}

It is easy to check that slope stability is invariant if $L$ is replaced by some power, and thus makes sense when $L$ is an ample $\mathbb Q$-divisor.  In many cases it is possible to calculate the slope of $X$ and $Z$
using the Grothendieck-Riemann-Roch formula, and for the purpose of this
paper we will mostly be concerned with the case that $Z$ is an effective
divisor in a smooth surface $X$ in which case 
\begin{align}\label{eq:slopeofdivisorinsurface}
  \mu(X,L) &= -\frac{K_X.L}{L^2}\quad\text{and}\\
  \mu_c(\O_Z,L) &= \frac{3(2L.Z-c(K_X.Z+Z^2))}{2c(3L.Z- cZ^2)}
\end{align}
where $K_X$ denotes the canonical class of $X$ \cite[Thm.\ 5.3]{ross_thomas:06:obstr_to_exist_const_scalar}.   Note that stability of $(X,L)$ with respect to an effective divisor $Z$ depends only on the numerical equivalence class of $L$ and $Z$.   We stress however that a destabilising divisor is necessarily integral.


\chapter{Exceptional Divisors}\label{sec:exceptional}

In this section we discuss how the existence of exceptional divisors of high genus
in a surface witnesses instability.  

\begin{defn} Let $D=\sum_{i=1}^m d_i
  D_i$ be an effective divisor in a smooth surface with irreducible components $D_i$ for $i=1,\dots,m$.  We say $D$ is \emph{exceptional} if the intersection matrix
  $(D_i.D_j)$ is negative definite.
\end{defn}

\begin{thm}\label{thm:negativedestabilises}
  Suppose that $X$ is a smooth projective surface containing an
  exceptional divisor $D$ that has arithmetic genus at least two and such that
  \begin{equation}
  D.D_i\le 0 \text{ for every irreducible component } D_i \text{ of }D .\label{eq:definitionofnegative}
\end{equation}
  Then $D$ destabilises $(X,L)$ for some polarisation $L$.
\end{thm}

\begin{proof}
  Let $D=\sum_i d_i D_i$ where $D_i$ are the irreducible components, and fix an ample divisor $H$.   Since the matrix $(D_i.D_j)$ is negative definite, there exist rational numbers $q_i$ such that
\[L_0 := H+ \sum_i q_i D_i\]
satisfies $L_0.D_j=0$ for all $j$.  If we write $\sum q_i D_i = E-F$ where $E$ and $F$ are effective with no common components then clearly $E.F-F^2\ge 0$ as $F^2\le 0$, whereas $E.F-F^2\le 0$ by definition of the $q_i$.  Thus $F=0$ so $\sum q_i D_i$ is effective and, since its intersection with each $D_j$ is strictly negative, we deduce that each $q_i$ is strictly positive. 
Let \[\varepsilon=\operatorname{min}_i \left\{\frac{q_i}{d_i}\right\}\] which is strictly positive.
Then $L_0-cD$ is nef for all $0\le c\le \varepsilon$, for
\eqref{eq:definitionofnegative} implies $(L_0-cD).D_j\ge 0$ for all $j$ and if $C$ is an
irreducible curve not supported on $D$ then $(L_0-cD).C=(H + \sum_i
(q_i-cd_i)D_i).C\ge 0$.  We note that $L_0-cD$ need not be
ample due to the possible existence of $D_j$ with $D.D_j=0$.

Now set
\[ L_s := L_0 + sH,\]
which is ample for positive $s$.  We will show that $D$ destabilises $(X,L_s)$ for sufficiently small $s$.  To this end note first that $L_s-cD$ is ample for $0<c\le \varepsilon$, which implies the Seshadri
constant satisfies $\varepsilon(D,L_s)\ge \varepsilon$.    Now $L_0^2= H.L_0\ge H^2>0$, so
$\mu(X,L_0)=-K_X.L_0/L_0^2$ is finite.  On the other hand as $L_0.D=0$,
\[\mu_c(\O_D,L_0)=\frac{3(2L_0.D-c(K_X.D+D^2))}{2c(3L_0.D-cD^2)}=\frac{3(2p_a(D)-2)}{2c
  D^2}.\] Hence $\mu_c(\O_D,L_0)$ tends to minus infinity as $c$ tends
to $0$.  Thus there is a $0<c<\varepsilon$ such that $ \mu_c(\O_D,L_0)<
\mu(X,L_0)$.  So as long as $s>0$ is sufficiently small we have
$c< \varepsilon(D,L_s)$ and $ \mu_c(\O_D,L_s)< \mu(X,L_s)$.  Hence
 $D$ destabilises $(X,L_s)$ as claimed.
\end{proof}

\begin{cor}\label{cor:exceptionaldestabilises}
  Suppose $X$ is a smooth projective surface containing an exceptional divisor $D$ that has arithmetic genus at least two.  Then $(X,L)$ is slope unstable for some polarisation $L$.    In particular this holds if $X$ contains an irreducible divisor with arithmetic genus at least two and negative self-intersection.
\end{cor}

The proof centres around the numerical cycle of $D$ whose definition we now recall.  Suppose $D=\sum d_i D_i$ is exceptional.    Then the fact that the intersection matrix is negative definite
  implies there exists a (non-trivial) effective divisor $D'
  =\sum d'_i D_i$ which satisfies $D'.D_i\le 0$ for
  all $i$.  Moreover it is easy to check that if $D'$ and $D''$ are two
  such divisors, then $\text{min}(D',D'')$ is another.
  Thus, as long as $D$ is connected,  there is a unique minimal divisor with this property which is called the \emph{numerical cycle} of $D$ and is denoted
  by $D_{\num}$ (see \cite{reid:97:chapt_algeb_surfac} for an introduction to this circle of ideas).

\begin{proof}[Proof of \linkref{Corollary}{cor:exceptionaldestabilises}]
Without loss of generality we may assume $D$ is connected.  Then Artin's well known classification of rational singularities states $p_a(D_{\num})\ge 0$, and if $p_a(D_{\num})=0$ then $p_a(D')\le 0$ for all effective divisors $D'$ supported on $D$ \cite[Thm.\ 3]{artin:66:isolat_ration_singul_surfac}.  

Similarly a result of Laufer  \cite[Cor.\ 4.2]{laufer:77:minim_ellip_singul} says that if $p_a(D_{\num})=1$ then $p_a(D')\le 1$ for all such $D'$ (see  also \cite{nemethi:99:ellip_goren_singul_surfac}).  Hence the assumptions imply that $p_a(D_{\num})\ge 2$, so the Corollary follows from \linkref{Theorem}{thm:negativedestabilises}. 
\end{proof}

\begin{rmk}\label{rmk:singular}
  The divisor $L_0$ used in the above proof satisfies $L_0.C=0$ if and
  only if $C$ is supported on $D$.  Suppose that there exists a contraction $X\to Y$ of $D$ (by Grauert's Theorem \cite{grauert:62:uber_modif_und_exzep_analy_mengen} this is always possible in the analytic category) and that $L_0$ is the pullback of a divisor $L_Y$.  For small $s$, the divisor $L_s$ makes the volume of $D$ 
  small, and the fact that $D$ destabilises $(X,L_s)$  can be
  interpreted as saying that $(Y,L_Y)$ is destabilised by the
  singularity obtained by contracting $D$.      The Corollary then says that for $Y$ to be stable it can have only rational and elliptic singularities.  This is analogous to \cite[Prop.\ 3.20]{mumford(77):stabil_projec_variet} which considers the case of asymptotic Chow stability and it is worth pointing out that similar results have been obtained by Ishii \cite{ishii:83:chow_instab_certain_projec_variet}, again in the case of Chow stability.  It would be interesting to have a list of the possible singularities on a K-stable surface.
\end{rmk}

\begin{example}
  There is a wealth of examples of surfaces that contain exceptional 
  curves of negative self-intersection and arithmetic genus at least
  two.  These include:
  \begin{enumerate}
  \item If $C$ is a smooth curve of genus at least 2 then the diagonal
    in $C\times C$ has negative self-intersection; moreover for
    certain non-general curves $C$ there is another irreducible curve of
    negative self-intersection in $C\times C$ that has arithmetic genus at least two.
    The instability of these examples are those of
    \cite{shu:unstab_kodair_fibrat} and
    \cite{ross:06:unstab_produc_smoot_curves} respectively.
  \item By blowing up sufficiently many points on an irreducible curve
    of genus at least 2 in a smooth surface one obtains a surface $X$
    with an irreducible curve of negative self intersection and genus at least 2.
    Moreover by taking a branched cover of $X$ one obtains in this way
    examples that are of general type.
  \end{enumerate}
  The reader will no doubt be able to provide many more examples.
\end{example}

\begin{example}
  The assumption in \linkref{Theorem}{thm:negativedestabilises} that the
  arithmetic genus be at least 2 is necessary.  For example, a K3
  surface is slope semistable with respect to every polarisation,  since it has trivial canonical
  class \cite[Thm.\ 8.4]{ross_thomas:07:study_hilber_mumfor_criter_for}, and this is true irrespective of
  the existence of any $-2$ curves.  The next result says that this is typical.
\end{example}

\begin{prop}\label{prop:genus}
  Suppose that $(X,L)$ is a smooth polarised surface with $K_X.L\ge 0$ (which holds for every $L$ if $X$ has non-negative Kodaira dimension).    If $D$ is an effective divisor that destabilises $(X,L)$ then $D$ has arithmetic genus at least 2.
\end{prop}
\begin{proof}
  By assumption,  $\mu(X,L)=-K_X.L/L^2\le 0$.  On the other hand if $D$ is an effective divisor then
  \[\mu_c(\O_D,L) =
  \frac{3(2L.D-c(K_X.D+D^2))}{2c(3L.D-cD^2)}>-\frac{3c(2p_a(D)-2)}{2c(3L.D-cD^2)}.\]
  Since the denominator is positive for $0<c\le \varepsilon(D,L)$ we conclude that $\mu(X,L)>
  \mu_c(\O_D,L)$ implies $p_a(D)\ge 2$.
\end{proof}

\begin{example}
  A surface with negative Kodaira dimension can be destabilised
  by a divisor with arithmetic genus less than 2.  For instance $\PP^2$
  blown up at a single point is destabilised by the exceptional
  divisor for all polarisations \cite[Ex.\ 5.27]{ross_thomas:06:obstr_to_exist_const_scalar}
\end{example}

  We now give an example of a destabilising divisor that is not exceptional.  However, in this as well as other examples we have studied, it is the case that if $D$ destabilises then there is an exceptional $D'\subset D$ that also destabilises.  It would be interesting to know if this is a general phenomena.

\begin{example}\label{ex:nonexcep}
  Let $C_g$ be a smooth curve of genus $g\ge 2$ and $C_h$ be a smooth curve of genus $h=(g-1)d+1$ that admits an unramified covering $p\colon C_h\to C_g$ of degree $d$.   Set $X = C_g\times C_h$, and let $E$ be the graph of $p$.  Furthermore let $F_g$ be a fibre that has genus $g$ and similarly for $F_h$.  The canonical class of $X$ is $K_X = (2g-2) (F_h + dF_g)$ and since $E$ has genus $h$ a simple calculation with the adjunction formula implies $E^2 = d(2-2g)$.  Let
\[ L = (2g-2)dF_g +E \quad\text{and}\quad D = F_g + E.\]
Clearly $D$ is effective, $L$ is nef and $L.E=0$.  Moreover $\varepsilon(D,L)\ge 1$.  An easy calculation yields 
\[\mu(X,L) = -2g \qquad \mu_c(\O_D,L) = \frac{3(2-c(d(2g-2)+2g))}{2c(3-c(d(2-2g)+2))}.\]
Thus $\mu_c(\O_D)\to -\frac{3}{2c}$ as $d\to \infty$.  Setting $c=3/(8g)<\varepsilon(D,L)$ this implies  that $D$ destabilises $(X,L')$ for $d\gg 0$, as long as $L'$ is a polarisation sufficiently close to $L$.  On the other hand it is clear that $D$ is not exceptional as $F_g^2=0$.  A similar calculation shows that $E$ also destabilises when $d$ is large.
\end{example}

\chapter{Reduction to Divisors}\label{sec:reduction}

In this section fix a smooth polarised surface $(X,L)$.  Our aim is the following theorem.

\begin{thm}\label{thm:reduction}
  If $(X,L)$ is slope unstable then it is destabilised by a divisor.
\end{thm}

To prove this we deform an arbitrary destabilising subscheme $Z$ into a subscheme $Z'$ which is a \emph{disjoint} union of a divisor and a zero-dimensional piece.  By semicontinuity we conclude that $Z'$ also destabilises, and thus some connected component of $Z'$ destabilises.  The proof is completed by showing that zero-dimensional subscheme can never destabilise a smooth surface.

\begin{prop}\label{prop:connectedcomponent}\cite{ross_thomas:07:study_hilber_mumfor_criter_for}
  Suppose that $Z$ destabilises a polarised variety $(X,L)$.  Then there is a connected component of $Z$ that also destabilises.  Moreover if a thickening $jZ$ destabilises for some $j\ge 1$ then so does $Z$.
\end{prop}
\begin{proof}
  This is taken from \cite[Prop.\ 4.25]{ross_thomas:07:study_hilber_mumfor_criter_for}.  Suppose $Z=Z_1\cup Z_2$ where $Z_1\cap Z_2$ is empty.  It is easy to check that $\varepsilon(Z)\le \varepsilon(Z_i)$ for $i=1,2$ (this can be done using the Kleiman criterion on the blowup of $Z$).   Moreover since $Z_1\cap Z_2$ is empty we can write $\mu_c(\O_{Z_i}) = \alpha_i/\beta_i$ where $\beta_i>0$ and $\mu_c(\O_{Z_1\cup Z_2}) = (\alpha_1+\alpha_2)/(\beta_1+\beta_2)$.  Thus if $Z$ destabilises then either $Z_1$ or $Z_2$ must destabilise.  The second statement follows from elementary properties of the slope formula, and the observation that $\epsilon(jZ) = j^{-1}\epsilon(Z)$.
\end{proof}

\begin{prop}\label{prop:disjoint}
Suppose that $(X,L)$ is slope unstable.  Then either it is destabilised by a divisor or by a zero-dimensional subscheme.
\end{prop}
\begin{proof}
Let $Z$ be a destabilising subscheme of $X$.  If $Z$ is zero dimensional then we are done.  Otherwise let $D$ be the largest divisor contained in $Z$ so we can write $\I_Z=\I_D\cdot \I_{Y_0}$ where $Y_0$ is a zero-dimensional subscheme.  Explicitly, if $Z$ is defined locally by an ideal $I$ and $D$ is defined by a function $f$ then $Y_0$ is defined locally by the ideal $(I:f) = \{ g: gf\in I\}$.

Clearly it is possible to move $Y_0$ away from its support in such a way that the thickenings of the family are all flat.  That is, there exists a small disc $0\in \Delta\subset \mathbb C$ and a subscheme $Y\subset X\times \Delta$ such the fibre of $Y$ over $0$ is $Y_0$, and such that the thickenings $jY$ are all flat over $\Delta$.  Moreover by choosing $Y$ generically and shrinking $\Delta$ if necessary we can assume that $Y_t$ is disjoint from $D$ for all $t\neq 0$. (This is easiest seen in the analytic topology where one can take general local analytic coordinates $x,y$ around each point of the support of $Y_0$ and move the component of $Y_0$ supported at this point along a coordinate axis by $(x,y)\mapsto (x+t,y)$ for $t\in\Delta$).

Now consider the subscheme $W$ given by $\I_W = \I_{D\times \Delta}\cdot \I_{Y}$ and let $\mathcal X\to X\times \Delta$ be the blowup along $W$ with exceptional divisor $\mathcal E$.  Since $\mathcal X$ is isomorphic to the blowup along $Y$ we have that $\mathcal X$ is flat over $\Delta$.  The central fibre of $\mathcal X$ is isomorphic to the blowup of $X$ along $Z$, and $\mathcal E$ restricts to the exceptional divisor of this blowup.  Similarly the fibre of $\mathcal X$ over $t\neq 0$ is isomorphic to the blowup of $X$ along $Y_t\cup D$ (and again the $\mathcal E$ restricts to the exceptional divisor).

Now set $Z'=W_t$ for some $t\neq 0$.  Then $\mu_c(\O_{Z}) = \mu_c(\O_{Z'})$ for all $c$ by flatness of $\mathcal X$.  Moreover by semicontinuity of Seshadri constants (since ampleness is an open condition), $\epsilon(Z)\le \epsilon(Z')$.  Thus if $Z$ destabilises $X$ then so does $Z'$.  Hence from \linkref{Proposition}{prop:connectedcomponent} some connected component of $Z'$ also destabilises and this component is either a divisor or zero-dimensional.
\end{proof}

We now turn to stability with respect to zero-dimensional subschemes.

\begin{lem}\label{lem:slopeofzerodim}
If $Z\subset X$ is supported at a point then
\begin{equation}
  \mu_c(\O_Z) \ge \frac{3}{2c}.\label{eq:slopesupportedpoint}
\end{equation}

\end{lem}
\begin{proof}
 Let $\pi\colon  \tilde{X}\to X$ be the blowup along $Z$ with exceptional divisor $E$ and  $p\colon \bar{X}\to \tilde{X}$ be a resolution of singularities.  We
  set $q=\pi \circ p\colon \bar{X}\to X$ and $F=p^* E$.  Note that if $0<x<\epsilon(Z)$ then $\pi^*L-xE$ is nef, and thus so is $q^*L-xF$.  Now applying the Grothendieck-Riemann-Roch formula to $\bar{X}$ and using \cite[Sec.\ 8.1]{ross_thomas:07:study_hilber_mumfor_criter_for} we have
\[ \mu_c(\O_D) = \frac{\int_0^c \tilde{a}_1(x) + \frac{\tilde{a}_0'(x)}{2} dx }{\int_0^c \tilde{a}_0(x)dx} \]
where for all $0<x<\varepsilon(Z)$,
\begin{align*}\label{eqn:boundsonai}
\tilde{a}_0(x) = -\frac{x^2 F^2}{2} \quad\text{and}\quad \tilde{a}_1(x)\ge -\frac{K_X.L}{2}+\frac{K_{\bar{X}}}{2} . (q^*L-xF)
\end{align*}
(to compare this with \cite[4.20, 8.2, 8.3]{ross_thomas:07:study_hilber_mumfor_criter_for} use  $a_i(x) = a_i -\tilde{a}_i(x)$).  We have $K_{\bar{X}}-q^* K_X\ge 0$ so $\tilde{a}_1(x)\ge 0$ whenever $q^*L-xF$ is nef, and thus $\tilde{a}_1(x)\ge 0$ for $0<x<c$.  Integration now yields the desired inequality.
\end{proof}

\begin{rmk}\
  \begin{enumerate}
  \item With a little more work it can be shown that the inequality \eqref{eq:slopesupportedpoint} is always strict.  In particular this is the case when $Z=\{p\}$ is a single point in which case $\mu_c(\O_Z) = 3/c$.  On the other hand the slope $\mu_c(\O_{mZ})$ of a $m$-thickened point gets arbitrarily close to $3/(2c)$ as $m$ tends to infinity.

  \item   At this point it is easy to prove  \linkref{Theorem}{thm:reduction} when $K_X.L\ge 0$ (which occurs, for instance, whenever $X$ has non-negative Kodaira dimension),  since then $\mu(X)\le 0$ whereas
  the slope of any zero-dimensional subscheme is strictly positive.
\end{enumerate}
\end{rmk}

The next two Lemmas give bounds on the Seshadri constant of zero-dimen\-sional subschemes.

\begin{lem}\label{lem:seshadrisupportedpoint}
Let $Z_1\subset Z_2$ be subschemes of $X$ supported at a point.  Then
\[\epsilon(Z_1)\ge \epsilon(Z_2).\]
\end{lem}
\begin{proof}
This lemma is presumably well known to experts, but for convenience we include a proof here.  
Recall that a sheaf $\I$ is said to be $m$-regular (with respect to $L$) if
\[H^i(\I\otimes L^{m-i})=0 \quad \text{for all }i>0,\]
and that the Mumford-Castelnuovo regularity is defined to be 
\[\reg(\I) = \operatorname{min}\{m : \I \text{ is $m$-regular}\}.\]
The relation with Seshadri constants is given by a theorem of Cutkosky-Ein-Lazarsfeld \cite[Thm.\ B]{cutkosky-ein-lazarsfeld(01):posit_compl_ideal_sheav} which says that if a subscheme $Z$ is defined by the ideal sheaf $\I_Z$ then
\begin{equation}
\frac{1}{\epsilon(Z)} = \lim_{j\to \infty} \frac{1}{j}\reg(\I_Z^j)\label{eq:cel}.
\end{equation}

We claim that if $0\to \I_2\to \I_1\to Q\to 0$ is an exact sequence of sheaves with $Q$ supported at a point then $\reg(\I_2)\ge \reg(\I_1)$.   But this is clear as $H^i(Q\otimes L^{m-i})=0$ for all $m$ and all $i>0$, so there is a surjection
\[H^i(\I_2\otimes L^{m-i})\to H^i(\I_1\otimes L^{m-i})\to 0.\]
Thus if $\I_2$ is $m$-regular then so is $\I_1$, and hence $\reg(\I_2)\ge \reg(\I_1)$.  Applying this to the inclusion $\I_{Z_2}^j\subset \I_{Z_1}^j$,  this implies $\reg(\I_{Z_2}^j) \ge \reg(\I_{Z_1}^j)$ for all $j\ge 1$, so the lemma follows from \eqref{eq:cel}. 
\end{proof}

\begin{rmk}
  It is important in the above lemma that $Z_1$ and $Z_2$ be supported at a point.  For example let $X$ be the blowup of $\PP^2$ at a point with exceptional divisor $E$.  Denote the hyperplane class by $H$ and let $L=3H-E$ which is ample.  Suppose $Z_1$ is the proper transform of a line through the point blown up, and $Z_2$ is the union of $Z_1$ with the exceptional divisor.  Then numerically $Z_1\equiv H-E$ and $Z_2\equiv H$ and one can check $\epsilon(Z_1,L)=1$ whereas $\epsilon(Z_2,L)=2$.
\end{rmk}

 The next result comes from \cite{ross_thomas:07:study_hilber_mumfor_criter_for} where it was used to show that reduced points do not destabilise.

\begin{lem}\label{lem:reducedpoint}
  Let $Z$ be a reduced point $Z=\{p\}$.  Then
\[ (3L + c K_X).L\ge 0 \quad \text{for all } 0<c\le \epsilon(Z).\]
\end{lem}
\begin{proof}
See the proof of \cite[Thm. 4.29]{ross_thomas:07:study_hilber_mumfor_criter_for}.
\end{proof}

We now complete the proof of the main theorem. 

\begin{proof}[Proof of \linkref{Theorem}{thm:reduction}]
We claim that no zero-dimensional subscheme $Z\subset X$ destabilises.  By \eqref{prop:connectedcomponent} we may assume without loss of generality that $Z$ is connected and thus supported at a point $p$. For $m\ge 1$ we write $\epsilon(mp)$ for the Seshadri constant of the $m$-thickening of the point $p$ (i.e. the subscheme given by $\I_{\{p\}}^m$).    Suppose first that $mp\subset Z$ for some $m\ge 2$.  Then by \eqref{lem:seshadrisupportedpoint},
\[\epsilon(Z)\le \epsilon(mp) = \frac{1}{m}\epsilon(p) \le \frac{1}{2}\epsilon(p).\]
If $c\le \epsilon(Z)$ then $2c\le \epsilon(p)$ so by \eqref{lem:reducedpoint},
\[ (3L + 2cK_X).L\ge 0.\]
Hence from \eqref{lem:slopeofzerodim},
\[\mu_c(\O_Z) \ge \frac{3}{2c} \ge \frac{-K_X.L}{L^2} = \mu(X)\]
and so $Z$ does not destabilise.  

Thus we may suppose that $Z$ does not contain $mp$ for any $m\ge 2$.      To calculate the slope of such a $Z$ note first that $\mu_c(\O_Z)$ depends only on the local geometry of $X$ around $p$ (in the analytic topology), so we may work in the ring $R=\mathbb C[[x,y]]$ and suppose that $Z$ is defined by an ideal $I\subset R$.    Since $I$ is not contained in $(x,y)^2$ there is a change of coordinates such that $I$ is of the form
\[ I = (x,y^n)\]
for some $n\ge 1$.   Then $R/I^j$ is spanned by $x^{\alpha}y^{\beta}$ for $0\le \alpha <j$, $0\le \beta< (j-\alpha)n$, so
\[\operatorname{length}(R/I^{j}) = \sum_{\alpha=0}^{j-1} n(j-\alpha)=\frac{nj(j+1)}{2}.\]
Replacing $j$ by $xk$ (where $xk\in \mathbb N$) and using the definition of $\tilde{a}_i(x)$ from Section \ref{sec:slopestability} this implies
\[ \tilde{a}_0(x) = \frac{nx^2}{2}  \quad\text{and}\quad    \tilde{a}_1(x) = \frac{nx}{2}\]
which yields
\[\mu_c(\O_Z)=\frac{3}{c}.\]

Now as $Z$ is supported at $p$, we have $\epsilon(Z)\le\epsilon(p)$ \eqref{lem:seshadrisupportedpoint}. Thus from \eqref{lem:reducedpoint}, if $0<c\le \epsilon(Z)$ then
\[\mu_c(\O_Z) = \frac{3}{c} \ge -\frac{K_X.L}{L^2} = \mu(X)\]
so $Z$ does not destabilise.  Hence no zero-dimensional subscheme can destabilise, and \linkref{Theorem}{thm:reduction} follows from \linkref{Proposition}{prop:disjoint}.
\end{proof}

\chapter{Destabilising divisors on minimal surfaces}\label{sec:nonnegativeI}

We now turn to the task of showing that any divisor that destabilises a surface with non-negative Kodaira dimension has negative self-inter\-section.  In this section we prove a theorem that deals with the case of minimal surfaces.

\begin{thm}\label{thm:negativedivisor}
Let $(X,L)$ be a smooth polarised surface with non-negative Kodaira dimension.  Suppose $D$ is a divisor in $X$ with $D^2\ge 0$ and that either $4K_X.D + K_X^2\ge 0$ or $K_X.D\le 0$.  Then $D$ does not destabilise $(X,L)$.
\end{thm}

\begin{cor}\label{cor:negativeminimal}
  Let $(X,L)$ be a polarised surface with non-negative Kodaira dimension such that $K_X^2\ge 0$ (in particular this holds for minimal surfaces).   Then any divisor that destabilises $(X,L)$ must have negative self-inter\-section.  
\end{cor}

The proof consists of an elementary analysis of the slope of a divisor with non-negative self-intersection.

\begin{defn} We say an effective divisor $D\subset X$ \emph{pseudo-destabilises} a polarised surface $(X,L)$ if 
\[\mu_c(\O_D)<\mu(X) \text{ for some } 0<c<\tilde{\varepsilon}(D,L)\]
where  $\tilde{\varepsilon}(D,L)$ is defined to be the smallest positive root of $p(t) = (L-tD)^2.$
\end{defn}

\begin{rmk}
(1) Clearly $\varepsilon(D,L)\le \tilde{\varepsilon}(D,L)$, so if $D$ destabilises then it pseudo-destabilises.  (2) Using the Hodge index theorem, one can check $\tilde{\varepsilon}(D,L)$ is well defined, and $D^2\tilde{\varepsilon}(D,L)\le L.D$.  Thus $\mu_c(\O_D)$ is finite for $0<c\le \tilde{\varepsilon}(D,L)$.  (3)  Whether a divisor pseudo-destabilises depends only on its numerical equivalence class.
\end{rmk}

\linkref{Theorem}{thm:negativedivisor} is a consequence of the following strengthening.

\begin{thm}\label{thm:negativedivisorII}
Let $(X,L)$ be a smooth polarised surface with non-negative Kodaira dimension.  Suppose $D$ is a divisor in $X$ with $D^2\ge 0$ and that either  $4K_X.D + K_X^2\ge 0$ or $K_X.D\le 0$.  Then $D$ does not pseudo-destabilise $(X,L)$.
\end{thm}
\begin{proof}
As $X$ has non-negative Kodaira dimension, $\mu(X)\le 0$. If $p_a(D)\le 1$ then $\mu_c(\O_D)>0$ for all $0<c<\tilde{\epsilon}(D)$ so $D$ does not pseudo-destabilise.  Thus we may assume $p_a(D)\ge 2$ at which point it is elementary to check that $\mu_c(\O_D)$ is strictly decreasing in the range $0<c<\tilde{\varepsilon}(D)$.  Thus to prove the theorem it is sufficient to show
\begin{equation}
\mu(X,L)\le \mu_{\tilde{\varepsilon}(D)}(\O_D,L). \label{eq:slopeinequalitytobeprovedII}
\end{equation}

Now since the statement of the theorem is invariant under rescaling $L$, we may assume without loss of generality that $L.D=L^2$.  If $D=L$ then $\tilde{\varepsilon}(D)=1$ and since $K_X.L\ge 0$,
\[ \mu_{\tilde{\varepsilon}(D)}(\O_D) = \frac{3(L^2-K_X.L)}{4L^2}> \frac{-K_X.L}{L^2} = \mu(X),\]
so $D$ does not pseudo-destabilise.  If $D\neq L$  then consider the plane in $N^1(X)$ spanned by $D$ and $L$.  Since the restriction of the intersection form to this plane has index $(1,-1)$ we can write
\[D = L + y \tau\]
where $\tau$ is a numerical class that is orthogonal to $L$ and satisfies $\tau^2=-L^2$.

 By assumption $0\le D^2 = L^2(1-y^2)$ so  $|y|\le 1$, and changing the sign of $\tau$ if necessary we may assume $0\le y\le 1$.  A simple calculation yields 
\[\tilde{\varepsilon}:=\tilde{\varepsilon}(D)=(1+y)^{-1}.\]

Let $K'=zL + w\tau$ be the projection of $K_X$ to the plane $\Lambda$.  Then $K_X.L = K'.L=zL^2$, so $z\ge 0$,  and $K_X.D = K'.D=L^2(z-wy)$.  Hence
\begin{align*}
\mu_{\tilde{\varepsilon}}(\O_D,L) &= \frac{3(2L.D-\tilde{\varepsilon}(K_X.D+D^2))}{2\tilde{\varepsilon}(3L.D-\tilde{\varepsilon} D^2)}= \frac{3(2-\tilde{\varepsilon}(z-yw+1-y^2))}{2\tilde{\varepsilon}(3-\tilde{\varepsilon}(1-y^2))}\\
&= \frac{3(2(1+y) - (z-yw) - 1+y^2)}{2(2+y)},\\
\mu(X,L) &= -\frac{K_X.L}{L^2}=-z.
\end{align*}
Now if $w\ge 0$ then,
\begin{align*}
  \mu_{\tilde{\varepsilon}}(\O_D,L) &= \frac{3(2(1+y) - (z-yw) - 1+y^2)}{2(2+y)}\ge \frac{3((1+y)^2 - z)}{2(2+y)}\\
&\ge -\frac{3z}{2(2+y)}\ge -z,
\end{align*}
so $D$ does not pseudo-destabilise. Next note that if $K_X.D\le 0$ then $w\ge 0$ so we are done by the above.   Thus we are left to deal with the case that $w<0$ and $4K_X.D + K_X^2\ge 0$ in which case
\begin{align*}
0&\le  4K_X.D + K_X^2 \le 4K'.D + K'^2 \\
&= L^2[4(z-wy)+ z^2-w^2]\le L^2[4(z-w)+z^2 - w^2 ] \\
&= L^2(z-w)(z+w+4)
\end{align*}
which implies $w\ge -z-4$.  Hence
\begin{align*}
\mu_{\tilde{\varepsilon}}(\O_D,L) &= \frac{3(2(1+y) - z+yw - 1+y^2)}{2(2+y)}
\ge\frac{3((1-y)^2- z(1+y))}{2(2+y)}\\
&\ge  -\frac{3z(1+y)}{2(2+y)}\ge -z = \mu(X,L).
\end{align*}
Thus $D$ does not pseudo-destabilise $(X,L)$ as claimed. 
\end{proof}

\begin{example}
Let $g_1,g_2\ge 2$ and suppose $C_1$ and $C_2$ are very general curves of genus $g_1$ and $g_2$ respectively  (by this we mean $(C_1,C_2)$ is a very general point in the product $M_{g_1}\times
  M_{g_2}$ of moduli spaces) and set $X=C_1\times C_2$.  Then clearly $K_X$ is ample, and it is well known that $N^1(X)$ has rank 2.  Thus $X$ has no curves of negative self-intersection so is slope stable with respect to all polarisations.
We remark that this result is expected from considerations of stability of products.  If $(X,L)$ and $(X',L')$ are K-stable, then one expects that $(X\times X',L\otimes L')$ is also stable (this is motivated by the conjectured equivalence between K-stability and the existence of constant scalar curvature K\"ahler metrics).    
\end{example}

\begin{example}\label{ex:weinkove}
Let $X$ be a surface with ample canonical class, and $L$ be a polarisation such that
\begin{equation}
2(K_X.L)L- (L^2)K_X \text{ is positive.}\label{eq:weinkovecondition}
\end{equation}
We claim that $(X,L)$ is slope-stable.  Since $K_X$ is ample we only have to check slope stability with respect to divisors $D$ with negative self-intersection. But this is clear, for if $0<c\le \varepsilon(D,L)$ then $2LD-cD^2> (L-cD).D\ge 0$ so
\begin{align*}
  \mu_c(\O_D,L) &= \frac{3(2L.D-cK_X.D-cD^2)}{2c(3LD-cD^2)} > -\frac{3K_X.D}{2(3L.D-cD^2)}\\
&\ge -\frac{K_X.D}{2L.D}\ge -\frac{K_X.L}{L^2}=\mu(X,L).
\end{align*}
\end{example}
\begin{rmk}
The condition \eqref{eq:weinkovecondition} was introduced by Chen \cite{chen(00):lower_bound_mabuc_energ_its_applic} and Donaldson  \cite{donaldson(99):momen_maps_diffeom} regarding the existence of special metrics on surfaces.  It was studied by Weinkove and Song-Weinkove \cite{weinkove:04:conver_j_flow_kaehl_surfac,song_weinkove:conver_singul_j_flow_with} who prove that this condition is sufficient for the Mabuchi functional to be proper, which is expected to imply K-stability.  \linkref{Example}{ex:weinkove} is an improvement to the surface case of \cite[Thm.\ 5.5]{ross_thomas:06:obstr_to_exist_const_scalar} that contains a similar result for higher dimensions.
\end{rmk}

\chapter{Destabilising divisors on surfaces with non-negative Kodaira dimension}\label{sec:nonnegativeII}

This section is devoted to proving the following:

\begin{thm}\label{thm:negativedoesnotdestab}
  Let $(X,L)$ be a smooth polarised surface with non-negative Kodaira dimension.  Then any divisor that destabilises $(X,L)$ has negative self-intersection.
\end{thm}

The proof proceeds by discarding certain components of $D$ in such a way that, after a possible blowdown of $X$, the divisor $K_X+D$ is nef at which point the results from the previous section are applied.

\begin{prop}\label{prop:adjoint}
  Let $(X,L)$ be a polarised surface of non-negative Kodaira dimension and suppose that a divisor $D$ pseudo-destabilises $(X,L)$. Then there exists a (possibly trivial) blowdown $\pi: X\to X'$, an effective divisor $D'$ and polarisation $L'$ of $X'$ such that $D'$ pseudo-destabilises $(X',L')$ and 

  \begin{enumerate}
  \item $D-\pi^*D'$ is a sum of components of $D$ that have negative self-intersection and $D'^2\ge D^2$,
  \item $L=\pi^* L'-E_L$ where $E_L$ is supported on the exceptional divisor of the blowdown \emph{and}
  \item $K_{X'}+D'$ is nef.   
  \end{enumerate}
\end{prop}

Granted this, it is easy to prove the main theorem of this section.

\begin{proof}[Proof of \linkref{Theorem}{thm:negativedoesnotdestab}]
  Suppose for contradiction that $D$ is a divisor with $D^2\ge 0$ that destabilises.  Applying the proposition we obtain a divisor $D'$ that pseudo-destabilises a surface $X'$ with $D'^2\ge 0$ and $K_{X'} +D'$ nef.  But $X'$ also has non-negative Kodaira dimension so $0\le (K_{X'}+D').K_{X'}$.  Thus
either $K_{X'}D'\le 0$ or $4K_{X'}.{D'} + K_{X'}^2\ge 0$ which is impossible by \linkref{Theorem}{thm:negativedivisorII}.
\end{proof}

Thus it remains to prove the Proposition which we start now by simplifying destabilising divisors.

\begin{lem}\label{lem:removenegative}
  Let $(X,L)$ be a smooth polarised surface with non-negative Kodaira dimension.  Suppose that $D$ is an effective divisor that is numerically equivalent to $D'+F$ where $D'$ and $F$ are non-trivial effective divisors with $F^2<0$ and
  \begin{equation}
2D.F \le \operatorname{min}\{F^2, F^2-K_X.F\}.\label{eq:conditionforremoval}
\end{equation}
\begin{enumerate}
\item If $D$ pseudo-destabilises then so does $D'$.
\item If in addition $D.F_i\le F_i^2$ for every irreducible component $F_i$ of $F$ and $D$ destabilises then so does $D'$.
\end{enumerate}
\end{lem}

  To understand this statement, let $D$ be a divisor that
  pseudo-destabilises and suppose $D$ contains an effective divisor
  $D'$ such that  $D'^2\ge D^2$ and $p_a(D')\ge p_a(D)$.
  Then the Lemma implies $D'$ also pseudo-destabilises provided $F=D-D'$ satisfies $F^2<0$.  Although this formulation is simpler, it is slightly weaker and is not the form in which
  \linkref{Lemma}{lem:removenegative} will be used. In fact we will
  be interested in the case that $D=D'+F$ where $F$ is an $-n$ curve (i.e.\ $F\simeq \PP^1$ with $F^2=-n$), in which case
  condition \eqref{eq:conditionforremoval} becomes
  \begin{align*}
    D.F \le F^2 +1 &\quad \text{if $F$ is a $-n$ curve with $n\ge 2$}\\
    D.F \le -1 &  \quad \text{if $F$ is a $-1$ curve.}
  \end{align*}

\begin{proof}[Proof of \linkref{Lemma}{lem:removenegative}]
Notice that $2D'.F + F^2= 2DF-F^2\le 0$, so for $c>0$
\[ (L-cD)^2 = (L-cD' -cF)^2 \le (L-cD')^2 + c^2 (2D'.F+F^2) \le (L-cD')^2\]
which implies $\tilde{\varepsilon}(D')\ge \tilde{\varepsilon}(D)$.  Now the assumption that $D$ pseudo-destabilises means $\mu_c(\O_D)<\mu(X)$ for some $0<c<\tilde{\varepsilon}(D)$.  But
\begin{align*}
  \mu_c(\O_D) &= \frac{3(2L.D -c (K_X.D+D^2))}{2c(3L.D-cD^2)} \\
&=\frac{3(2L.D' -c (K_X.D'+D'^2)) + 3(2L.F -c(K_X.F+F^2+ 2D'.F))}{2c(3L.D'-cD'^2) + 2c(3L.F-cF^2-2cD'.F)}
\end{align*}
and by hypothesis
\[\frac{3(2L.F -c(K_X.F+F^2+ 2D'.F))}{2c(3L.F-cF^2-2cD'.F)} > 0\ge \mu(X).\]
Recall now that if $\alpha,\beta,\mu\in \mathbb R$ and $\gamma,\delta\in\mathbb R^+$ with $(\alpha+\beta)/(\gamma+\delta)<\mu$ and $\beta/\delta>\mu$ then $\alpha/\gamma<\mu$.  Thus
\[\mu_c(\O_{D'})=\frac{3(2L.D' -c (K_X.D'+D'^2))}{2c(3L.D'-cD'^2)}<\mu(X)\]
so $D'$ pseudo-destabilises $(X,L)$ as claimed.

For the second statement it is sufficient to prove that $\varepsilon(D)\le \varepsilon(D')$.  But if $L-cD=L-cD'-cF$ is ample, then clearly $(L-cD').C>0$ for any irreducible curve $C$ not contained in $F$, and if $F_i$ is a component of $F$ then $(L-cD').F_i \ge L.F_i - c(D.F_i - F_i^2)>0$, so $L-cD'$ is also ample.
\end{proof}

\begin{proof}[Proof of \linkref{Proposition}{prop:adjoint}]  We claim that, without loss of generality, we may suppose that if $F$ is an irreducible curve with $F^2<0$ then
  \begin{equation}\label{eq:canbeassumed}
\begin{split}
D.F\ge 0  &\quad\text{ if $F$ is a $-1$ curve, and } \\
(K_X+D).F\ge 0 &\quad\text{ otherwise. }
\end{split}
\end{equation}
To see this suppose $D$ pseudo-destabilises and \eqref{eq:canbeassumed} does not hold for some irreducible $F$ with $F^2<0$.   If $F$ is not a $-1$ curve then $D.F< -K_X.F\le 0$, and if $F$ is a $-1$ curve then $D.F<0$ by hypothesis, so either way $F$ must be a component of $D$.   Suppose first $p_a(F)\ge 1$.  Then $(K_X+D).F<0$ and $K_X.F\ge 0$ so
\[ 2D.F = 2(K_{X}+D).F - (2p_a(F)-2) + F^2 -K_X.F\le F^2 - K_{X}.F.\]
 Thus from \linkref{Lemma}{lem:removenegative} we deduce that $D'=D-F$ also pseudo-destabilises.  Similarly if $D$ is a $-n$ curve with $n\ge 2$ then $D.F\le -K_X.F -1 = F^2+1$, and if $F$ is a $-1$ curve then $D.F\le -1$, so once again the Lemma implies $D'=D-F$ pseudo-destabilises.    One can easily check that in each of these cases $D'^2\ge D^2$, and thus by making $D$ smaller if necessary  we  can assume \eqref{eq:canbeassumed}.

Suppose now that $F$ is a $-1$ curve with $D.F=0$.  Let $\pi\colon X\to Y$ be the blowdown of $F$, so $D\equiv \pi^* D_Y$ for some effective divisor $D_Y$ in $Y$.  Moreover $L=\pi^* L_Y -sF$ for some polarisation $L_Y$ of $Y$ and $s>0$, which implies $\mu(X,L)<\mu(Y,L_Y)$.  Now $D.F=0$ implies $\mu_c(\O_{D},L)=\mu_c(\O_{D_Y},L_Y)$.   Moreover if $c>0$ then $(L-cD)^2\le (L_Y-cD_Y)^2$ so $\tilde{\varepsilon}(D_Y,L_Y) \ge \tilde{\varepsilon}(D,L)$.    Thus $D_Y$ pseudo-destabilises $Y$ and $D_Y^2=D^2\ge 0$.    Thus by repeating this process we obtain a divisor $D'$ that pseudo-destabilises a polarised surface $(X',L')$ such that $(K_{X'}+D').F\ge 0$ for all irreducible curves $F$ with negative self-intersection.  But if $F$ is irreducible with non-negative self-intersection then it is nef, so $(K_{X'}+D').F\ge 0$ since $K_{X'}$ and $D'$ are effective.  Thus $K_{X'}+D'$ is nef as claimed.
\end{proof}

\chapter{Stability with respect to nef divisors}\label{sec:nef}

The purpose of this section is to prove that nef divisors do not destabilise.     Of course when  $X$ has non-negative Kodaira dimension this has already been proved in the previous section.   However the proof we give here is different, and applies as well to surfaces of negative Kodaira dimension.

\begin{thm}\label{thm:nef}
Let $X$ be a smooth polarised surface and $D$ be a nef divisor.  Then $D$ does not destabilise $X$ with respect to any polarisation.
\end{thm}

The rough idea for the proof is that, apart from two special cases, there exists a blowdown $\pi\colon X\to X'$ such that $D$ is the pullback of a divisor $D'$ and $K_{X'}+2D'$ is nef.  (The two special can only arise when $X$ is a ruled surface, or a blowup of $\PP^2$.)   We then prove slope stability of $X$ by comparison with $X'$.   


\begin{prop}\label{3cases}
Let $X$ be a smooth surface and $D$ be a nef divisor. There exists 
a (possibly trivial) blowdown $\pi: X\to X'$ to a smooth surface $X'$ and a divisor $D'\subset X'$
such that $D=\pi^*D'$ and one of the following conditions hold:
\begin{enumerate}
\item The divisor $K_{X'}+2D'$ is nef.
\item  $X'=\mathbb P^2$ and $D'$ is a line.
\item  The surface $X'$ is a $\mathbb P^1$-bundle over a curve and $D'$ is a union of its
  fibres.
\end{enumerate}
Moreover in the first two cases the divisor $K_X+3D$ is pseudoeffective (i.e.\ it lies in the closure of the effective cone).
\end{prop}

\begin{proof}
We will show that under the assumption that $D$ 
intersects positively all $-1$ rational curves in $X$ the pair
$(X,D)$ is already of these three types.   The proposition clearly follows from this 
because we can subsequently blowdown all  $-1$ curves  in $X$
that have zero intersection with $D$. 

 So suppose that $D.C\ge 1$ for all $-1$ curves $C$ and that $K_X+2D$ is not nef, so there exists an irreducible curve $C$ with $(K_X+2D).C<0$, which implies $K_X.C<0$. By the Mori Cone theorem \cite{mori:82:threef_whose_canon_bundl_are} $C=\sum_i R_i$,  
where each $R_i$ is an extremal ray of $X$. Since $(K_X+2D).C<0$ at least one of the rays $R_i$ is not a $-1$ curve. This 
is possible only if $X$ is a $\mathbb P^1$-bundle over a curve, or $X$ is $\mathbb P^2$.

If $X$ is a $\mathbb P^1$-bundle over a curve then $D$ must be a union of its 
fibres otherwise $K_X+2D$ is nef. 
Indeed all $\mathbb P^1$-bundles over curves apart from $\mathbb P^1\times \mathbb P^1$
and $\mathbb P^2$ blown up at one point have just 
one extremal ray $R$ namely the fibre of the surface. 
If $D$ is not a fibre then $D.R>0$ and so $(K_X+2D).R=-2+2D.R\ge 0$. 
The cases of $\mathbb P^1\times \mathbb P^1$ and $\mathbb P^2$
blown up at one point are similar.
If $X$ is $\mathbb P^2$ then clearly $D$ is a line, otherwise $K_X+2D$ is nef.

For the last statement notice that in the first two cases the divisor $K_{X'}+3D'$ is nef, and hence so is $\pi^*(K_{X'}+3D')$.  Since $K_{X}-\pi^* K_{X'}$ is effective, this implies $K_X+3D$ is pseudoeffective.
\end{proof}

\begin{rmk}
  Of course Cases 2 and 3 can only occur when $X$ has negative Kodaira dimension, but we stress that the converse is not true.  Also it is clear that if $X$ has non-negative Kodaira dimension then the proof actually produces a pair $(X',D')$  such that $K_{X'}+D'$ is nef. 
\end{rmk}

\begin{proof}[Proof of \linkref{Theorem}{thm:nef}]

Fix a polarisation $L$ on $X$ and a nef divisor $D$.   The aim is to show that $D$ does not destabilise, which entails showing that if $c$ is chosen so $A=L-cD$ is ample then
\[\mu(X,L)=-\frac{K_X.L}{L^2} <
\frac{3(2L.D-c(K_X.D+D^2))}{2c(3L.D-cD^2)}=\mu_c(\O_D).\]
Rearranging, the goal becomes to prove that the polynomial
\begin{align*} \label{eq:definitionofP}
P(c)=&c^3(K_X.D+3D^2)D^2+4c^2(K_X.A+3A.D)D^2 \nonumber\\
&+ 3c[2((K_X+2D).A)(A.D)-(K_X.D-D^2)A^2]+6(A.D)A^2
\end{align*}
is positive.  Since $A$ is ample, the constant term is clearly positive, so is sufficient to prove that all the other coefficients of $P(c)$ are non-negative.  At this point it is convenient to define
\begin{equation}
 Q(A,D) = 2((K_X+2D).A)(A.D)-(K_X.D-D^2)A^2\label{eq:definitionofQ}.
\end{equation}

 Let $\pi\colon X\to X'$ be the blowdown coming from \linkref{Proposition}{3cases}, so $D=\pi^* D'$ for some divisor $D'$.   Notice that in the first two cases of \eqref{3cases}, $K_X+3D$ is pseudoeffective, and since $D$ and $A$ are nef this implies the highest two coefficients of $P(c)$ are non-negative.   In the third case, $D^2=0$ and so these coefficients vanish.  Thus is sufficient to prove that $Q(A,D)\ge 0$ which we do for each case separately. \smallspace

\noindent \emph{Case 1:} (The divisor $K_{X'}+2D'$ is nef)   

Note that if $K_X.D-D^2\le 0$ then clearly $Q(A,D)\ge 0$, so we may suppose $K_X.D-D^2>0$.   Write $A=\pi^*A'-E_L$, where $E_L$ is supported on the exceptional divisor and $A'$ is an ample divisor on $X'$.  Then
\begin{align*}
Q(A,D)&=2((K_X+2D).A)(A.D)-(K_X.D-D^2)A^2\\
&=Q(A',D')+2((K_X-\pi^*K_{X'}).A)(A.D)-(K_X.D-D^2)E_L^2.
\end{align*}
Obviously $E_L^2\le 0$ and since $K_X-\pi^*K_{X'}$
is effective we have $(K_X-\pi^*K_{X'}).A\ge 0$.   Thus $Q(A',D')\le Q(A,D)$, so the proof of this case is completed by the following lemma.

\begin{lem}\label{lem:Q}
   Let $X'$ be a smooth surface.  Suppose that $D'$ is a pseudo-effective divisor with $D'^2\ge 0$ such that $K_{X'}+2D'$ is pseudoeffective and $(K_{X'}+2D')^2\ge 0$.    Then $Q(A',D')\ge 0$ for every ample divisor $A'$.
\end{lem}

\begin{proof}

We will prove the result for any effective $D'\in N^1(X)\otimes \mathbb R$.  As above we clearly can assume $K_{X'}.D'-D'^2>0$.  Suppose first that the numerical class of $D'$ is proportional to $K_{X'}+2D'$.  Then $(K_{X'}+2D')=\alpha D'$ for some $\alpha\ge 0$ and $Q(A',D') = 2\alpha (A'.D')^2 - (\alpha-3)A'^2D'^2$ which is non-negative by the Hodge-index theorem.  

Hence we can assume $K_{X'}+2D'$ and $D'$ span a 2-plane $\Lambda\subset N^1(X)\otimes \mathbb R$.  Since both divisors are effective and have non-negative square,  the
restriction of the intersection form to $\Lambda$ has signature  $(1,-1)$
and moreover the vectors $D'$ and $K_{X'}+2D'$ are in the same non-negative octant.

Now notice that if $\tilde{A}$ is the projection of $A'$ onto $\Lambda$ then $\tilde{A}^2\ge A'^2$ so $Q(\tilde{A},D)\le Q(A',D)$.    Thus we may assume without loss of generality that $A'$ lies in $\Lambda$.     Moreover by homogeneity of $Q(A',D)$ with respect to $A'$ we can assume $A'^2=1$.  Thus by picking a suitable basis for $\Lambda$ we can write $D'=(x,-cx)$ and $K_{X'}+2D'=(z,cz)$  with $x,z\ge 0$, $0\le c\le 1$ and $A'=(\sqrt{1+y^2},y)$.  Hence
 \begin{align*}
Q(A',D')&=2((K_{X'}+2D').A')(A'.D')-(K_{X'}.D'-D'^2)A'^2\\
&=2((K_{X'}+2D').A')(A'.D')-((K_{X'}+2D').D'-3D'^2)A'^2\\
&=2z(\sqrt{1+y^2}-cy)x(\sqrt{1+y^2}+cy)-(zx+c^2zx-3(x^2-c^2x^2))\\
&\ge2zx(1+y^2-c^2y^2)-(zx+c^2zx)
\ge zx-c^2zx\\
&\ge 0.
\end{align*}
as required.
\end{proof}

To deal with the remaining two cases we need an additional lemma.

\begin{lem}\label{lem:exceptionalinequality} 
Let $p\colon X\to Y$ be a blowdown of smooth surfaces whose exceptional divisor $E$ is connected.  Fix a smooth divisor $D_Y$ in $Y$ such that $p(E)\in D_Y$ and let $A = p^* A_Y - E_A$ be an ample divisor on $X$ with $E_A$ supported on $E$.   Finally let $\tilde{D}_Y$ be the proper transform of $D_Y$ and set
\[ F=p^* D_Y - \tilde{D}_Y.\]
Then
\[(A.K_{p})(A.F)\ge -E_A^2\]
where $K_{p}=K_X-p^* K_Y$ is the relative canonical class of $p$.
\end{lem}

\begin{proof}
We argue by induction on the number of irreducible components of $E$.  If $E$ is irreducible it must be a $-1$ curve, in which case it is clear that equality holds.  Suppose now that $E$ has more than one irreducible component, so at least one component $E_0$ is a $-1$ curve.  Let $q: X\to Z$ be the contraction of $E_0$, and $p_Z\colon Z\to Y$ be the remaining blowdown, so $p=p_Z\circ q$.   We have $A=q^* A_Z - aE_0$ for some ample divisor $A_Z$ on $Z$.     

Let $F_Z$ be the difference between $p_Z^* D_Y$ and the proper transform of $D_Y$ in $Z$.  We claim that
\begin{equation}\label{equation:ainequality}
A.F \ge A_Z. F_{Z}>a.
\end{equation}
To see this let $x=q(E_0)$.  If $x$ is not contained in the proper transform of $D_Y$ we have $F=q^*(F_Z)$, and otherwise $q^*(F_Z)$ is strictly contained
in $F$ which proves the first inequality.  For the second inequality let $G$
be any irreducible component of $F_Z$ that contains $x$, and $\tilde{G}$ be its proper transform in $X$.   Since $A$ and $A_Z$ are ample and $\tilde{G}.E_0\ge 1$,
\[0<A.\tilde{G}=(q^*A_Z-aE_0).\tilde{G}\le A_Z.G-a\le A_Z.F_Z-a.\]

Now we finish the proof of the lemma. Writing $A_Z = p_Z^* A_Y- E_{A_Z}$ 
where $E_{A_Z}$ is supported on the exceptional divisor of $p_Z$,  we have by induction 
$(A_Z.K_{p_Z})(A_Z.F_Z) \ge - E_{A_Z}^2$.  The identity $K_{p}=q^*K_{p_Z}+E_0$ and \eqref{equation:ainequality} thus  imply
\begin{align*}
 (A.K_{p})(A.F)  &\ge (A_Z.K_{p_Z} + a)(A_Z.F_Z) \ge 
  -E_{A_Z}^2 + aA_Z.F_Z \\
&\ge -E_{A_Z}^2 + a^2= -E_A^2
\end{align*}
as required.
\end{proof}

\noindent We now conclude the proof of \linkref{Theorem}{thm:nef}.\smallspace

\noindent \emph{Case 2: } ($X'=\mathbb P^2$ and $D'$ is a line)

We will show $Q(A,D)> 0$ by induction on the number of \emph{connected} components $m$ of the exceptional divisor $E$ of the blowdown $X\to X'$.   If $m=0$ then $X=\PP^2$ and $D=D'$ is a line at which point $Q(A,D) = 2A^2>0$.   Suppose now that $m=1$ and $A=\pi^* A'-E_A$ where $E_A$ is supported on $E$.  Since $Q(A,D)$ depends only on the numerical equivalence class of $D$ we may without loss of generality assume that $D=\pi^* D'$ where $D'$ is a line that contains the point $\pi(E)$.  Then
\begin{align*}
 Q(A,D) &= Q(A',D') + 2(A. (K_X-\pi^* K_{X'})) (A.D) + 4 E_A^2 \\
&\ge Q(A',D')+ 2E_A^2  
&\tag*{by \linkref{Lemma}{lem:exceptionalinequality}} \\
&= 2A'^2 + 2E_A^2 \tag*{as $Q(A',D') = 2A'^2$}\\
&=2 A^2>0.
\end{align*}

For the inductive step let $E_1$ and $E_2$ be two connected components of $E$ and suppose $x_1=\pi(E_1)$ and $x_2=\pi(E_2)$ are their images in $\PP^2$.  Let $D'$ be the line through $x_1$ and $x_2$, and $\tilde{D'}$ be its proper transform in $X$.  Then $\pi^*D'-\tilde{D'}$ is a collection of exceptional divisors, and two of them $F_1$ and $F_2$ are
supported on $E_1$ and $E_2$ respectively. Without loss of 
generality we can assume that $A.F_1\le A.F_2$, and since $D\equiv \pi^*D' $ and  $\pi^*D'-F_1-F_2$ is effective, this implies $A.D\ge 2A.F_1$. 

Now let $p\colon X\to Y$ be the blowdown of $E_1$ 
and  $D_Y$ be the pullback of $D'$ under the remaining blowdown
$Y\to X'$.  We write $A=p^* A_Y - E_A$ where $E_A$ 
is supported on $E_1$.  Noticing that in this case $K_{X'}D'-D'^2=-4$,
\begin{align*}
Q(A,D) &= Q(A_Y, D_Y)  + 2(A.(K_X-p^* K_Y)) (A.D) + 4 E_A^2\\
 &\ge Q(A_Y, D_Y)  + 4(A.(K_X-p^* K_Y))(A.F_1) + 4 E_A^2\\
&\ge  Q(A_Y,D_Y)  \tag*{ by \linkref{Lemma}{lem:exceptionalinequality}}\\
&>0 \qquad \tag*{ by inductive hypothesis}
\end{align*}
which completes the proof of Case 2. \smallspace

\noindent \emph{Case 3:} ($X'$ is a $\mathbb P^1$-bundle over a curve and $D'$ is a union of fibres)

If $D_1$ and $D_2$ are effective divisors with $D_1.D_2=0$ then  $Q(A,D_1+D_2)\ge Q(A,D_1) + Q(A,D_2)$.  Thus as all fibres of $X'$ are numerically equivalent we may assume that $D'$ is a single fibre.  Notice that in this case $K_{X}.D-D^2=-2$.

Suppose that $F$ is a connected component of $E$, and let $p\colon X\to Y$ be the blowdown along $F$.    Write $A=p^* A_Y - F_A$ and let $D_Y$ be the pullback of $D'$ to $Y$.    Again since $Q(A,D)$ depends on the numerical equivalence class of $D$ we may assume that $D_Y$ contains $p(F)$.  Thus by \linkref{Lemma}{lem:exceptionalinequality} we get
  \begin{align*}
Q(A,D)&=Q(A_Y,D_Y)+2(A.(K_X-p^*K_{Y}))(A.D)+2 F_{A}^2 \ge Q(A,D_Y).  
\end{align*}
Thus by induction on the number of connected components of $E$ we can suppose that $X$ is $\mathbb P^1$-bundle over a curve and $D$ is a fibre, in which case a simple calculation yields $Q(A,D)\ge 0$.  This finishes the proof of Case 3, and completes the proof of \linkref{Theorem}{thm:nef}.
\end{proof}

\begin{example}\label{ex:p2blown2points}
Let $\pi\colon X\to \PP^2$ be the blowup of $\PP^2$ at two points.  We shall show that if $L=-K_X$ is the anticanonical polarisation then $(X,L)$ is slope stable but not K-stable.  In fact the latter statement is well known,  and follows from the existence of a $\mathbb C^*$-action which has non-trivial Futaki invariant (see for example \cite[Ex.\ 3.11]{tian00:canonical_kaehler}). 

To prove slope stability it is sufficient by \linkref{Theorem}{thm:reduction} to show $X$ is not destabilised by any divisor.    One way to do this is to start by checking directly that $X$ is not destabilised by any $-1$ curve (and thus by \eqref{prop:connectedcomponent} the same is true for any thickened $-1$ curve).  Assume, for contradiction, that a divisor $D$ destabilises,  so by \linkref{Theorem}{thm:nef} we know $D$ is not nef.  Thus there is a $-1$ curve $E$ in $X$ with $D.E<0$.  Let $p\colon X\to Y$ be the blowdown along $E$ (so either $Y$ is $\PP^1\times \PP^1$ or $Y$ is the blowup of $\PP^2$ at a single point).  Thus we can write $D=p^* D_Y + nE$ where $D_Y$ is an effective divisor in $Y$ and $n>0$.  

Now writing $L-cp^* D_Y = L-cD + ncE$ and using the Kleiman criterion we get that $\varepsilon(p^*D_Y)\ge \varepsilon(D)$.  Similarly writing $L-cnE = L-cD + cp^* D_Y$ we get $\varepsilon(nE)\ge \varepsilon(D)$ (this is clear when $Y=\PP^1\times \PP^1$ for then $D_Y$ is nef;  when $Y$ is $\PP^2$ blown up at a single point this follows as any exceptional component of $p^*D_Y$ is disjoint from $E$).  Thus as $p^* D_Y.E=0$ we deduce either $nE$ or $p^* D_Y$ destabilises.  But, as already observed, $nE$ does not destabilise and hence $p^* D_Y$ destabilises $X$.

Now $\varepsilon(D_Y,-K_Y)\ge \varepsilon(p^*D_Y,-K_X)$ whereas $\mu(X,-K_X)=\mu(Y,-K_Y)$ and $\mu_c(\O_{D_Y},-K_Y)=\mu_c(\O_{p^*D_Y},-K_X)$.  Thus $D_Y$ destabilises $(Y,-K_Y)$.  Now \linkref{Theorem}{thm:nef} implies this is impossible when $Y=\PP^1\times \PP^1$ for then $D_Y$ is nef.  On the other hand a direct calculation shows that the only destabilising divisor of $\PP^2$ blown up at a single point is the exceptional divisor.    Thus $D$ must be supported on the exceptional set of $\pi\colon X\to \PP^2$.  But from \eqref{prop:connectedcomponent}, if such a $D$ destabilises then so does one component of this exceptional set, which is not the case.  Thus $X$ is slope stable as claimed.

\end{example} 
\begin{rmk}
  Sz\'ekelyhidi has previously observed that, using the description of
  the slope of a toric subscheme from
  \cite[4.3]{ross_thomas:06:obstr_to_exist_const_scalar}, one can show
  that no toric subscheme destabilises $\PP^2$ blown up at two points.
  Since any effective divisor is numerically equivalent to a toric
  divisor, this also proves slope stability with respect to divisors.
\end{rmk}\smallspace\smallspace

\begin{flushright}
  {\small \scshape
Dmitri Panov, Dept.\ of Mathematics, Imperial College, \\ London, SW7 2AZ, UK \bigskip\\
Julius Ross, Dept.\ of Mathematics, Columbia University, \\New York, NY 10027, USA}
\end{flushright}
\clearpage

\bibliography{biblio2}
\end{document}